\documentclass[12pt,reqno,a4paper]{amsart}
\usepackage{amssymb,amscd, hyperref, epsf}

\def\hepsffile{\leavevmode\epsffile}

\begin{document}

\let\kappa=\varkappa
\let\epsilon=\varepsilon
\let\phi=\varphi
\let\p\partial

\def\Z{\mathbb Z}
\def\R{\mathbb R}
\def\C{\mathbb C}
\def\Q{\mathbb Q}
\def\P{\mathbb P}
\def\HH{\mathrm{H}}
\def\ss{{X}}

\def\conj{\overline}
\def\Beta{\mathrm{B}}
\def\const{\mathrm{const}}
\def\ov{\overline}
\def\wt{\widetilde}
\def\wh{\widehat}

\renewcommand{\Im}{\mathop{\mathrm{Im}}\nolimits}
\renewcommand{\Re}{\mathop{\mathrm{Re}}\nolimits}
\newcommand{\codim}{\mathop{\mathrm{codim}}\nolimits}
\newcommand{\id}{\mathop{\mathrm{id}}\nolimits}
\newcommand{\Aut}{\mathop{\mathrm{Aut}}\nolimits}
\newcommand{\lk}{\mathop{\mathrm{lk}}\nolimits}
\newcommand{\sign}{\mathop{\mathrm{sign}}\nolimits}
\newcommand{\rk}{\mathop{\mathrm{rk}}\nolimits}
\def\Jet{{\mathcal J}}

\renewcommand{\mod}{\mathrel{\mathrm{mod}}}

\newtheorem{mainthm}{Theorem}
\renewcommand{\themainthm}{{\Alph{mainthm}}}
\newtheorem{thm}{Theorem}[subsection]
\newtheorem{lem}[thm]{Lemma}
\newtheorem{prop}[thm]{Proposition}
\newtheorem{cor}[thm]{Corollary}

\theoremstyle{definition}
\newtheorem{exm}[thm]{Example}
\newtheorem{rem}[thm]{Remark}
\newtheorem{df}[thm]{Definition}

\renewcommand{\thesubsection}{\arabic{subsection}}
\numberwithin{equation}{subsection}

\title{Legendrian links, causality, and the Low conjecture}
\author[Chernov \& Nemirovski]{Vladimir Chernov and Stefan Nemirovski}
\thanks{The second author was supported by grants from DFG, RFBR, Russian Science Support Foundation,
and the programme ``Leading Scientific Schools of Russia.''}
\address{Department of Mathematics, 6188 Kemeny Hall,
Dartmouth College, Hanover, NH 03755, USA}
\email{Vladimir.Chernov@dartmouth.edu}
\address{%
Steklov Mathematical Institute, 119991 Moscow, Russia;\hfill\break
\strut\hspace{8 true pt} Mathematisches Institut, Ruhr-Universit\"at Bochum, 44780 Bochum, Germany}
\email{stefan@mi.ras.ru}

\begin{abstract} Let $(\ss^{m+1}, g)$ be a globally hyperbolic spacetime
with Cauchy surface diffeomorphic to an open subset of $\R^m$.
The Legendrian Low conjecture formulated by Nat\'ario and Tod says that
two events $x,y\in\ss$ are causally related if and only if the Legendrian
link of spheres $\mathfrak S_x, \mathfrak S_y$ whose points are light
geodesics passing through $x$ and $y$ is non-trivial in the contact manifold
of all light geodesics in $\ss$. The Low conjecture says that for $m=2$ the events
$x,y$ are causally related if and only if $\mathfrak S_x, \mathfrak S_y$
is non-trivial as a topological link. We prove the Low and the Legendrian Low conjectures. We also show
that similar statements hold for any
globally hyperbolic $(\ss^{m+1}, g)$ such that a cover 
of its Cauchy surface is diffeomorphic to an open domain in $\R^m.$
\end{abstract}

\maketitle

\vspace{-3pt}

\subsection{Introduction}
The space $\mathfrak N$ of non-parameterised future pointing null geodesics in a
globally hyperbolic spacetime $(\ss^{m+1},g)$, $m\ge 2$, has a natural structure
of a contact $(2m-1)$-manifold obtained by identifying $\mathfrak N$
with the spherical cotangent bundle $ST^*M$ of a smooth spacelike Cauchy surface $M^m\subset X$.
Null geodesics passing through a point~$x\in X$ form a Legendrian $(m-1)$-sphere
$\mathfrak S_x\subset\mathfrak N$ called the {\it sky\/} of~$x$.
(Details and definitions may be found in~\S\S\ref{lorentz}--\ref{nullgeo} below.)

All skies in~$\mathfrak N$ are Legendrian isotopic.
The situation is more interesting for {\it links\/} formed by pairs
of disjoint skies. It was observed by Low~\cite{Low0} that the isotopy class
of the link $\mathfrak S_x\sqcup\mathfrak S_y$ may depend on whether
the points $x$ and $y$ are {\it causally related\/},
that is, connected by a future pointing non-spacelike curve in~$\ss$.

It is not hard to show that all links $\mathfrak S_x\sqcup\mathfrak S_y$
formed by skies of causally {\it un\/}related points
belong to the same Legendrian isotopy class in~$\mathfrak N$ represented
by a pair of fibres of $ST^*M$.
It is therefore natural to call $\mathfrak S_x$ and $\mathfrak S_y$
{\it topologically unlinked\/} (respectively, {\it Legendrian unlinked\/})
if they are disjoint and the link $\mathfrak S_x\sqcup\mathfrak S_y$
is smoothly (respectively, Legendrian) isotopic to a link in that
`trivial' isotopy class.
It is also natural to ask whether the skies of causally related points
are in some sense linked. This question was raised in different forms
by Low~\cite{Low0}, \cite{Low1}, \cite{Low3}, \cite{LowLegendrian},~\cite{LowNullgeodesics}
and Nat\'ario and Tod~\cite{NatarioTod}. It appeared on Arnold's problem lists
as a problem communicated by Penrose~\cite[Problem~8]{ArnoldProblem},
\cite[Problem 1998-21]{ArnoldProblemBook}.
The following result was conjectured in~\cite[Conjecture 6.4]{NatarioTod} for the case when
the Cauchy surface is diffeomorphic to an open subset of $\R^3$.

\begin{mainthm}[Legendrian Low Conjecture]
\label{thma}
Assume that a smooth spacelike Cauchy surface of a globally hyperbolic spacetime $(X,g)$
has a cover diffeomorphic to an open subset of\/ $\R^m$, $m\ge 2$.
Then the skies of causally related points in $X$ are Legendrian linked.
\end{mainthm}

Explicit examples (see~\cite[\S 6]{NatarioTod}) show that
causally related points in a globally hyperbolic spacetime with Cauchy
surface diffeomorphic to $\R^3$ can have {\it topologically\/} unlinked skies.
That is, Legendrian linking is a strictly weaker condition for $m\ge 3$.
On the other hand, combining Theorem~\ref{thma} with a recent result
of Ding and Geiges~\cite{DG} on the classification of Legendrian links
in $ST^*\R^2$, we obtain the following result for $(2+1)$-dimensional
spacetimes. This result was conjectured by Low~\cite{Low0} for the case when
the Cauchy surface is diffeomorphic to an open subset of $\R^2$.

\begin{mainthm}[Low's Conjecture]
\label{thmb}
Assume that the universal cover of a smooth spacelike Cauchy surface of a globally hyperbolic
$(2+1)$-dimensional spacetime $(X,g)$ is diffeomorphic to\/ $\R^2$.
Then the skies of causally related points in $X$ are topologically linked.
\end{mainthm}

The proof of Theorem~\ref{thma} is based on the methods of the theory
of generating functions developed in the context of contact topology
by Traynor~\cite{Tr} and Bhupal~\cite{Bh} following the seminal work
of Viterbo~\cite{Vi}. (After our paper was submitted for publication,
we discovered that results very similar to those contained
in \S\S\ref{contact}--\ref{hodograph} below were obtained earlier
by Colin, Ferrand, and Pushkar'~\cite{CFP}.)
One more application of this approach shows that Legendrian linking,
unlike topological linking, can distinguish between past and future.
For the standard Minkowski spacetime, this result is an interpretation
of the main result of~\cite{Tr} in terms of skies, see \cite[Theorem 6.2]{NatarioTod}.

\begin{mainthm}
\label{thmc}
Assume that a smooth spacelike Cauchy surface of a globally hyperbolic spacetime $(X,g)$
has a cover diffeomorphic to an open subset of\/ $\R^m$, $m\ge 2$.
Let $x,y\in\ss$ be causally related points with disjoint skies.
Then the links $\mathfrak S_x\sqcup\mathfrak S_y$
and $\mathfrak S_y\sqcup\mathfrak S_x$ are not Legendrian isotopic.
\end{mainthm}

The universal cover of a $2$-dimensional manifold $M\neq S^2, \R\mathrm P^2$ with $\partial M=\varnothing$
is diffeomorphic to $\R^2$.
In particular, we see that Theorems~\ref{thma},~\ref{thmb},~\ref{thmc} hold
for $(2+1)$-dimensional globally hyperbolic spacetimes with Cauchy surfaces
other than $S^2$ or $\R\mathrm P^2$.

According to Agol's talk~\cite{AgolTalk}, Thurston's geometrisation
conjecture proved recently by Perelman~\cite{Perelman1},
\cite{Perelman2} implies that the universal cover of a closed
$3$-manifold $M$ is diffeomorphic either to $S^3$ or to an open
subset of $\R^3$. If $M$ is universally covered by $S^3$, then
the geometrisation conjecture tells us that $M$ is diffeomorphic
to a quotient of $S^3$ by a finite group of isometries of the standard round metric.
If $M$ is the interior of a compact $3$-manifold
$\ov M$ with boundary, then the double of $\ov M$ is a closed $3$-manifold
and hence $M$ is covered by an open subset of $\R^3$. Since
any open $3$-manifold can be exhausted by compact $3$-manifolds
with boundary, it follows (see Remarks~\ref{extended} and~\ref{exhaust})
that in the physically interesting case of $(3+1)$-dimensional globally
hyperbolic spacetimes, Theorems~\ref{thma} and~\ref{thmc}
hold assuming that the Cauchy surface is not diffeomorphic to
a metric quotient of the standard sphere~$S^3$.

If $M$ is a metric quotient of the standard round sphere $S^m$ by a finite group of isometries,
then it is easy to construct a globally hyperbolic spacetime with Cauchy surface~$M$
for which Theorems~\ref{thma}, \ref{thmb},~\ref{thmc} are false.
One can take the product Lorentz manifold $(M\times \R, \ov g\oplus -dt^2)$,
where~$\ov g$ is the quotient Riemann metric on~$M$, see~\cite[Example 3]{ChernovRudyak}.

It is worth pointing out that all known examples of this sort are {\it refocussing},
see~\cite{LowRefocussing} and~\cite[Definition 22]{ChernovRudyak}.
(This notion seems to be related to Riemann $Y_\ell^x$-manifolds studied
by B\'erard-Bergery~\cite{BerardBergery} and~Besse~\cite{Besse},
see~\cite[Remark 7]{ChernovRudyak}.)
On the other hand, it was proved by Rudyak and the first author~\cite[Corollary~1]{ChernovRudyak}
that if a globally hyperbolic spacetime is non-refocussing, then skies
of causally related points cannot be unlinked by a Legendrian isotopy
consisting of skies of points.
So it is conceivable that our results remain valid for any globally hyperbolic spacetime
that is not diffeomorphic to a refocussing spacetime.

\smallskip
\noindent
{\bf Contents of the paper.}
The key notion of non-negative Legendrian isotopy is introduced in \S\ref{legisotop}.
Necessary facts and definitions from Lorentz geometry are recalled in \S\ref{lorentz}.
Contact geometry of the space of null geodesics and its relation with
causality are discussed in~\S\ref{nullgeo}.
Generating functions are used to study Legendrian isotopies
in $1$-jet bundles in~\S\ref{contact}. The hodograph transformation
is applied in~\S\ref{hodograph}.
The proofs of the results stated in the introduction occupy \S\S\ref{proofleg}--\ref{prooftop}.
The last \S\ref{scifi} contains a tentative application to physics.

\smallskip
\noindent
{\bf Conventions.} All manifolds, maps etc.~are assumed to be smooth
unless the opposite is explicitly stated, and the word {\it smooth\/}
means $C^{\infty}$. The connected components of a disconnected manifold
(such as a link) are assumed to be ordered; maps between
disconnected manifolds are assumed to preserve the order of components.
Contactomorphisms of co-oriented contact structures are assumed
to be co-orientation preserving.

\subsection{Non-negative Legendrian isotopies}
\label{legisotop}
Let $Y$ be a contact manifold with a co-oriented
contact structure defined by a contact form~$\alpha$.
A submanifold $\Lambda\subset Y$ is called Legendrian if it is tangent to the contact distribution,
i.\,e., if $\alpha|_\Lambda\equiv 0$.
A {\it Legendrian isotopy\/} in $Y$ is a smooth family $\{\Lambda_t\}_{t\in[0,1]}$
of Legendrian submanifolds.
Two Legendrian submanifolds are called Legendrian isotopic if they
can be connected by a Legendrian isotopy.

A basic fact about Legendrian isotopies is the {\it Legendrian
isotopy extension theorem\/} (see, e.\,g., \cite[Theorem~2.6.2]{Ge}).
It asserts that for any Legendrian isotopy $\{\Lambda_t\}_{t\in[0,1]}$
of compact submanifolds, there exists a smooth family of
compactly supported contactomorphisms $\Psi_{t\in[0,1]}:Y\to Y$
such that $\Psi_0=\id_Y$ and $\Psi_t(\Lambda_0)=\Lambda_t$ for all~$t\in[0,1]$.
In particular, isotopic compact Legendrian submanifolds are isotopic
by an ambient contact isotopy.

\begin{df}
\label{defnonneg}
A Legendrian isotopy $\{\Lambda_t\}_{t\in[0,1]}$ in a contact manifold $(Y,\ker\alpha)$
is called {\it non-negative\/} if it has a parameterisation $F:\Lambda_0\times[0,1]\to Y$ 
such that $(F^*\alpha)\left(\frac{\p}{\p t}\right)\ge 0$.
\end{df}

Clearly, this definition does not depend on the choice of the parameterisation $F$
of the Legendrian isotopy and on the choice of the contact form defining the
co-oriented contact structure. It is also obvious that if $\Psi:Y\to Y'$
is a contactomorphism, then the image of a non-negative Legendrian isotopy in $Y$
is a non-negative Legendrian isotopy in~$Y'$.

\begin{exm}[Non-negative isotopy in $ST^*M$]
\label{isotopyremark}
Let $\pi:ST^*M\to M$ be the spherical cotangent bundle of an $m$-dimensional manifold $M$.
A point $p\in ST^*M$ may be regarded as a linear form $\wt p$ on $T_{\pi(p)}M$ defined
up to a multiplication by a positive scalar.
Thus, $p\in ST^*M$ is determined by the co-oriented hyperplane $\ell_p=\ker \wt p\subset T_{\pi (p)}M$.
(The co-orientation is given by the half-space in $T_{\pi (p)}M$ on which $\wt p$ is positive.)
The standard contact structure on $ST^*M$ is the co-oriented hyperplane distribution
$$
\chi=\{(d \pi)^{-1}(\ell_p)\subset T_p(ST^*M)\mid p\in ST^*M\}.
$$

Let $\{\Lambda_t\}_{t\in[0,1]}$ be a Legendrian isotopy parameterised by $F:\Lambda_0\times[0,1]\to ST^*M$.
It follows from the definitions that this isotopy is non-negative if and only if
the value of a linear form corresponding to $F(s, \tau)\in ST_{\pi(F(s,\tau))}^*M$
on the vector $d\pi\circ dF|_{(s,\tau)}(\frac{\partial}{\partial t})$ is non-negative
for all $(s,\tau)\in\Lambda_0\times[0,1]$.

Suppose now that the projection $\pi|_{\Lambda_t}:\Lambda_t\to M$ is an embedding for all $t\in[0,1]$.
For any point $p\in\Lambda_t$, the tangent hyperplane $T_{\pi(p)}\pi(\Lambda_t)$ coincides with $\ker\wt p$
and therefore has a canonical co-orientation. The co-oriented hypersurface $\pi(\Lambda_t)\subset M$
is called the wave front of the Legendrian submanifold $\Lambda_t\subset ST^*M$.
In this situation, the isotopy $\{\Lambda_t\}_{t\in[0,1]}$ is non-negative if and only if
the wave fronts $\pi(\Lambda_t)$ move in the direction of their co-orientation.\qed
\end{exm}

\begin{rem}
A closely related notion of non-negative contact isotopy was studied by
Eliashberg and Polterovich~\cite{ElPo}, Eliashberg, Kim, and Polterovich~\cite{ElKiPo},
and Bhupal~\cite{Bh}.
\end{rem}

\subsection{Lorentz geometry: definitions and terminology}
\label{lorentz}
Let $(\ss, g)$ be a Lorentz manifold of dimension $m+1$ and signature $(+,\dots,+,-)$.
A non-zero vector ${\bf v}\in T_p\ss$ is called {\it timelike, non-spacelike, null {\rm (}lightlike\/}), or
{\it spacelike\/}
if $g({\bf v}, {\bf v})$ is respectively negative, non-positive, zero, or positive.
A piecewise smooth curve is timelike if all of its velocity vectors are timelike.
Non-spacelike and null (lightlike) curves are defined similarly. Since $(\ss, g)$ has a unique
Levi-Civita connection, see for example~\cite[p.\,22]{BeemEhrlichEasley}, we can talk
about spacelike, timelike, and null (light) geodesics. A submanifold $M\subset \ss$ is
{\it spacelike\/} if $g$ restricted to $TM$ is a Riemann metric.

All non-spacelike vectors in $T_p\ss$ form a cone consisting of two hemicones, and
a continuous with respect to $p\in \ss$ choice of one of the two hemicones
is called a {\it time orientation\/} of $(\ss, g)$. The vectors from the chosen hemicones
are called {\it future pointing}. A time oriented connected Lorentz manifold is called a
{\it spacetime\/} and its points are called {\it events}.

For an event $x$ in a spacetime $(\ss, g)$, its {\it causal future\/} $J^+(x)\subset \ss$
(respectively, {\it chronological future\/} $I^+(x)$) is the set of all $y\in \ss$ that
can be reached by a future pointing non-spacelike (respectively, timelike) curve from~$x$.
The causal past $J^-(x)$ and the chronological past $I^-(x)$ of the event $x\in \ss$
are defined similarly.

Two events  $x,y$ are said to be {\it causally related\/} if $x\in J^+(y)$ or $y\in J^+(x)$;
and they are said to be {\it chronologically related\/} if and only if $x\in I^+(y)$
or $y\in I^+(x)$.

An open subset $U\subset \ss$ is itself a Lorentz manifold and for
$x\in U$ we denote by $J^{\pm}(x,U),$ $I^{\pm}(x,U)$ the causal and
the chronological past and future of $x$ with respect to~$U$.

\begin{df}
A spacetime is {\it causal\/} if it does not contain closed future pointing curves.
A spacetime $X$ is said to be {\it globally hyperbolic\/} if it is causal and
$J^+(x)\cap J^-(y)$ is compact for all $x,y\in \ss$.
\end{df}

This definition of global hyperbolicity is equivalent to the more classical one
used in~\cite{BeemEhrlichEasley} and~\cite{HawkingEllis}
by a recent result of Bernal and Sanchez~\cite[Theorem 3.2]{BernalSanchezCausal}.

A {\it Cauchy surface\/} in $(\ss, g)$ is a subset such that every inextensible future pointing
non-spacelike curve $\gamma(t)$ intersects it at exactly one value of $t.$
It is a classical result that $(\ss, g)$ is globally hyperbolic if and only if
it contains a Cauchy surface, see~\cite[pp.\,211--212]{HawkingEllis}.
Geroch~\cite{Geroch} proved that every globally hyperbolic $(\ss, g)$
is homeomorphic to the product of $\R$ and a Cauchy surface.
Bernal and Sanchez~\cite[Theorem 1]{BernalSanchez},%
~\cite[Theorem 1.1]{BernalSanchezMetricSplitting},~\cite[Theorem 1.2]{BernalSanchezFurther}
proved that every globally hyperbolic $(\ss^{m+1}, g)$ has a smooth spacelike Cauchy surface~$M^m$,
any two smooth spacelike Cauchy surfaces of $(\ss^{m+1},g)$ are diffeomorphic,
and that moreover for every smooth spacelike Cauchy surface $M$ there is a diffeomorphism
$h_M:M\times \R\to \ss$ such that
\begin{itemize}
\item[a)] $h_M(M\times t)$ is a smooth spacelike Cauchy surface for all $t$;
\item[b)] $h_M(x\times \R)$ is a future pointing timelike curve for all $x\in M$;
\item[c)] $h_M(M\times 0)=M$ with $h_M|_{M\times 0}:M\to M$ being the identity map.
\end{itemize}
This deep result has the following useful corollary
(cf.~\cite[Proof of Theorem 8]{ChernovRudyak}).

\begin{lem}
\label{move}
Let $x_1,x_2$ be causally unrelated points in~$\ss$.
Then they can be smoothly moved into $M$ so that
they remain causally unrelated in the process.
\end{lem}

\begin{proof}
Let $(\ov x_i,t_i)=h_M^{-1}(x_i)$ and assume that the points are ordered
so that $t_1\le t_2$. Note that moving $x_1$ into the future along the
segment $h_M(\ov x_1,[t_1,t_2])$ will not create causal relations
between $x_1$ and $x_2$. Indeed, if $h_M(\ov x_1,t)$ were causally
related to $x_2$ for some $t\in[t_1,t_2]$, then $x_2$ would have to
lie in $J^+(h_M(\ov x_1,t))$.
The union of a future pointing curve connecting $h_M(\ov x_1,t)$ to $x_2$
with the future pointing timelike curve $h_M(\ov x_1,[t_1,t])$ would
be a future pointing non-spacelike curve connecting $x_1$ to $x_2$.

Thus, we can move $x_1$ into the future till $t=t_2$ so that both points
lie on the Cauchy surface $h_M(M\times t_2)$. Since
points lying on the same Cauchy surface are causally unrelated,
it remains to use the obvious isotopy of Cauchy surfaces
$M_t:=h_M(M\times tt_2)$ connecting $h_M(M\times t_2)$ and
$M=h_M(M\times 0)$.
\end{proof}

\subsection{Contact geometry of the space of null geodesics}
\label{nullgeo}
Let $\mathfrak N$ be the space of all (inextensible) future pointing null geodesics
in a globally hyperbolic spacetime
$(\ss, g)$ considered up to affine orientation preserving reparameterisation.

\begin{df}
\label{skydef}
The {\it sky\/} $\mathfrak S_x$ is the set of all
null geodesics in $\mathfrak N$ passing through~$x\in\ss$.
\end{df}

Let $M$ be a spacelike Cauchy surface in $\ss$. A null geodesic $\gamma=\gamma(t)$ intersects $M$
at a time~$\ov t$.
Since $M$ is spacelike and $\gamma$ is a null geodesic, the linear form on $T_{\gamma(\ov t)}M$
given by $\mathbf v\mapsto g(\gamma'(\ov t),\mathbf v)$ is non-zero and therefore defines
a point in the spherical cotangent bundle $ST^*M$ of $M$. Thus, we have an identification
$$
\rho_M:\mathfrak N\overset{\simeq}{\longrightarrow} ST^*M.
$$
Note that if $\mathfrak S_x$ is the sky of a point $x\in M$, then
$\rho_M(\mathfrak S_x)=ST^*_xM$ is the fibre of $ST^*M$ over~$x$.

Consider the contact structure on $\mathfrak N$ induced by $\rho_M$ from the
standard contact structure on~$ST^*M$. Low~\cite{LowLegendrian} showed that
if $M'$ is another spacelike Cauchy surface, then the map $\rho_M\circ\rho_{M'}^{-1}:ST^*M'\to ST^*M$
is a {\it contactomorphism} (see also ~\cite[pp.\,252--253]{NatarioTod}).
Thus, this contact structure on $\mathfrak N$ does not depend on the choice
of the Cauchy surface~$M$.

It follows, in particular, that any sky $\mathfrak S_x$ is a {\it Legendrian\/}
sphere in $\mathfrak N$ because it corresponds to the fibre $ST_x^*M'$ for
a Cauchy surface $M'$ passing through~$x$. Note also that if two points
$x_1,x_2\in\ss$ are connected by a curve $x=x(t)$, then $\mathfrak S_{x(t)}$
is a Legendrian isotopy connecting the skies $\mathfrak S_{x_1}$ and $\mathfrak S_{x_2}$,
and hence all skies are Legendrian isotopic.

\begin{prop}
\label{causal}
Let $x,y\in\ss$ be causally related points not lying on the same null geodesic
and such that $y\in J^+(x)$. Then $\mathfrak S_y$ is connected to $\mathfrak S_x$
by a {\bf\em non-negative} Legendrian isotopy.
\end{prop}

\noindent
{\em Warning.}
Note that the non-negative isotopy goes to the past. This is caused by the
use of future pointing geodesics in the definition of $\rho_M$.

\begin{proof}
Since $x,y$ are not on the same null geodesic, we have $y\in I^+(x)$
by~\cite[Corollary 4.14]{BeemEhrlichEasley}. (It was Corollary~3.14 in the
first edition of~\cite{BeemEhrlichEasley}.)
Choose a past directed smooth timelike curve $\gamma:[0,1]\to \ss$ from $y$ to $x$.
We claim that the Legendrian isotopy $\{\mathfrak S_{\gamma(t)}\}_{t\in[0,1]}$
connecting $\mathfrak S_y$ to $\mathfrak S_x$ is non-negative.

Fix $\tau\in [0,1]$ and put $q=\gamma(\tau)$. Choose a spacelike Cauchy surface $M'$
of the type $h_{M}(M\times s)$ passing through~$q$. We will show that there exists
$\epsilon>0$ such that the Legendrian isotopy $\{\rho_{M'}(\mathfrak S_{\gamma(t)})\}$ in $ST^*M'$
is non-negative for all $t\in (\tau,\tau+\epsilon)$. The claim will follow because
$\rho_{M'}:\mathfrak N\to ST^*M'$ is a contactomorphism and $\tau$ can be
chosen arbitrarily.

\begin{figure}[htbp]
\begin{center}
\epsfxsize 12cm
\hepsffile{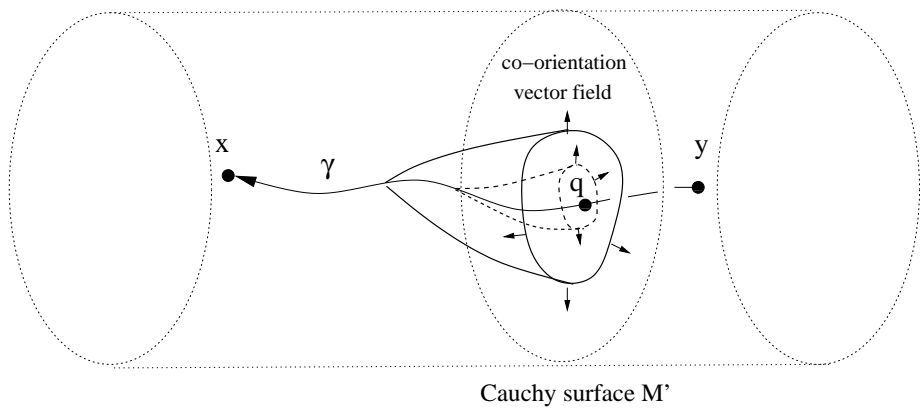}
\end{center}
\caption{Foliation of a neighbourhood of $q$ in $M'$ by embedded spheres}\label{ConeFoliation.fig}
\end{figure}

Using~\cite[Chapter 5, Proposition 7]{ONeill}, choose a convex neighbourhood $U$ of $q$ in $\ss$.
For every $p\in U$, the intersection of the exponent of the future null cone in $T_p\ss$
with $M'$ is either an embedded hypersurface or a single point in $M'$.
The intersection of $M'$ with the exponent of the null cone at $q$ is just one point $q$.
Moving $p$ from $q$ slightly into the past along $\gamma$, we see that the intersection of $M'$
with the exponent of the future null cone at $p$ is an embedded sphere $\Sigma_p\subset M'$.
Two such spheres corresponding to different points $q_1, q_2\in \gamma\cap U$ do not intersect.
Indeed, assume they do intersect at $z\in M'$. Without loss of generality assume
that $q_2\in I^+(q_1)$ and hence that $q_2\in I^+(q_1, U)$.
Then $z$ is on an arc of a null geodesic passing through~$q_1$.
The arc is contained in $U$ and thus $z$ is in $J^+(q_1, U)\setminus I^+(q_1, U)$
by~\cite[Chapter 14, Lemma 2]{ONeill}.
Similarly, $z$ is in $J^+(q_2, U)\setminus I^+(q_2, U)$ and, since $q_2\in I^+(q_1, U)$,
we have  $z\in I^+(q_1, U)$ by~\cite[Chapter 14, Corollary 1]{ONeill}. Contradiction.
Thus, a small deleted neighbourhood of $q$ in $M'$ is foliated by the spheres
$\Sigma_{\gamma(t)}$, $t\in (\tau,\tau+\epsilon)$,
that are expanding as the point $\gamma(t)$ moves into the past along~$\gamma,$
see Figure~\ref{ConeFoliation.fig}.

It follows from the definition of $\rho_{M'}$ that for any $t\in(\tau,\tau+\epsilon)$
the sphere $\Sigma_{\gamma(t)}\subset M'$ with the outward
co-orientation coincides with the wave front of the Legendrian submanifold
$\rho_{M'}(\mathfrak S_{\gamma(t)})\subset ST^*M'$.
Since the family of spheres is expanding, the wave fronts move in the direction of their
co-orientation and hence the Legendrian isotopy
$\{\rho_{M'}(\mathfrak S_{\gamma(t)})\}_{t\in (\tau,\tau+\epsilon)}$
is non-negative, see Example~\ref{isotopyremark}.
\end{proof}

The notion of (un)linking in $\mathfrak N$ is based on the following
observation going back to Low~\cite{Low0} (see also~\cite{NatarioTod}
and~\cite{ChernovRudyak}).

\begin{lem}
\label{unlinkclass}
The Legendrian isotopy class of the link $\mathfrak S_x\sqcup\mathfrak S_y\subset\mathfrak N$
formed by the skies of two causally unrelated points $x,y\in\ss$ does not depend on the points
$x$ and~$y$.
\end{lem}

\begin{proof}
Pick a spacelike Cauchy surface $M\subset\ss$. By Lemma~\ref{move}, we can move
$x$ and $y$ into $M$ keeping them causally unrelated. Since the skies
of causally unrelated points are disjoint, we obtain a Legendrian isotopy
connecting $\mathfrak S_x\sqcup\mathfrak S_y$ to a link of the form
$\rho_M^{-1}(ST^*_{x'} M\sqcup ST^*_{y'} M)$ for $x'\ne y'\in M$.
Once $M$ is fixed, any two such links are obviously Legendrian isotopic.
\end{proof}

\begin{df}
\label{unlinkdef}
Let $\mathcal U$ be the Legendrian isotopy class of links in $\mathfrak N$
containing $\mathfrak S_x\sqcup\mathfrak S_y$ for all causally unrelated points $x,y\in\ss$.
Two skies $\mathfrak S_x$ and $\mathfrak S_y$ are called {\it Legendrian unlinked\/} if they are disjoint
and $\mathfrak S_x\sqcup\mathfrak S_y$ belongs to~$\mathcal U$.
Two skies $\mathfrak S_x$ and $\mathfrak S_y$ are called {\it topologically unlinked\/}
if they are disjoint and $\mathfrak S_x\sqcup\mathfrak S_y$ is smoothly isotopic to a link in~$\mathcal U$.
Skies are called topologically (respectively, Legendrian) linked if they are not topologically
(respectively, Legendrian) unlinked.
\end{df}

Lemma~\ref{unlinkclass} shows that the Legendrian isotopy class $\mathcal U$ is well-defined.
A natural way to represent it in $\mathfrak N$ is to
identify $\mathfrak N$ with $ST^*M$ for a spacelike Cauchy surface $M\subset\ss$
and consider the link in $ST^*M$ formed by any two different fibres $ST_x^*M$
and $ST^*_yM$. In other words, we take the skies of two different points
$x\ne y$ lying on the same Cauchy surface.

\subsection{Generating functions and Legendrian isotopies}
\label{contact}
Let $\Jet^1(L)$ denote the $1$-jet bundle of a compact connected
manifold~$L$ equipped with the standard contact form $du-p\,dq$,
where $u$ is the fibre coordinate in $\Jet^0(L)$ and $p\,dq$
denotes the Liouville form on $T^*L$.

Let $\Lambda\subset\Jet^1(L)$ be a Legendrian submanifold.
A function $S=S(q,\xi):L\times\R^N\to\R$ is  a {\it generating function\/}
for $\Lambda$ if zero is a regular value of the partial differential $d_\xi S$
and the map
\begin{equation}
\label{gendef}
\{d_\xi S(q,\xi)=0\}\ni (q,\xi) \longmapsto (q,d_q S(q,\xi),S(q,\xi))\in \Jet^1(L)
\end{equation}
is a diffeomorphism onto~$\Lambda$. A generating function is said to be
{\it quadratic at infinity\/} if $S(q,\xi)=Q(q,\xi)+\sigma(q,\xi)$, where
$\sigma$ has compact support and $Q(q,\cdot)$ is a non-degenerate quadratic form in the variable~$\xi$.

Given a quadratic at infinity function $S:L\times\R^N\to\R$,
there is a topological procedure for selecting one of its critical values $c_-(S)$
(see~\cite[\S2]{Vi}). Consider the sublevel sets
$$
S^c:=\{(q,\xi)\in L\times\R^N\mid S(q,\xi)\le c\}
$$
and denote by $S^{-\infty}$ the set $S^c$ for a sufficiently negative $c\ll 0$.
Let $\kappa(S)$ denote the negative index of the quadratic form $Q(q,\cdot)$.
Pick a point $q_0\in L$ and a negative linear subspace $V\subset \{q_0\}\times \R^N$
for $Q(q_0,\cdot)$ of maximal possible dimension~$\kappa(S)$. The relative homology
class $[V]\in\HH_\kappa(L\times\R^N,S^{-\infty})$ does not depend on the choices made.
Define
$$
c_-(S):=\inf\{c\in\R\mid [V]\in \imath_*\HH_\kappa(S^c,S^{-\infty})\},
$$
where $\imath_*:\HH_\kappa(S^c,S^{-\infty})\to \HH_\kappa(L\times\R^N,S^{-\infty})$
is the homomorphism of relative homology groups
induced by the inclusion $\imath:S^c\to L\times\R^N$.

\begin{exm}
\label{graphs}
Let $\Lambda^f:=\{(q,df(q),f(q))\mid q\in L\}\subset\Jet^1(L)$ be the graph
of the $1$-jet of a smooth function $f:L\to\R$. It follows from the definitions
and Morse theory that
$$
c_-(S)=\min_{\{d_\xi S=0\}} S =\min_L f
$$
for any quadratic at infinity generating function $S:L\times\R^N\to\R$ of the Legendrian
submanifold~$\Lambda^f\subset\Jet^1(L)$. In particular, $c_-(S)=0$ for any quadratic at infinity generating
function of the zero section $\Lambda^0\subset\Jet^1(L)$.\qed
\end{exm}

Suppose now that $\{\Lambda_t\}_{t\in [0,1]}$ is a Legendrian isotopy such that $\Lambda_0$
is the zero section of $\Jet^1(L)$. By Chekanov's theorem~\cite{Che} (see also \cite{Cha} and~\cite{Tr})
there exists a smooth family of quadratic at infinity generating functions $S_t:L\times\R^N\to\R$ for $\Lambda_t$.

\begin{lem}
\label{nondecreasing}
If the Legendrian isotopy $\{\Lambda_t\}_{t\in [0,1]}$ is non-negative,
then $t\mapsto c_-(S_t)$ is a non-decreasing function on~$[0,1]$.
\end{lem}

\begin{proof}
It follows from the definition of a non-negative isotopy and formula~(\ref{gendef})
that $\frac{\p S_t}{\p t}(q,\xi)\ge 0$ when $d_\xi S_t(q,\xi)=0$ and,
in particular, when $d_{(q,\xi)}S_t(q,\xi)=0$. Hence, $c_-(S_t)$ is a
non-decreasing (continuous) function of $t$ by~\cite[Lemma~4.7]{Vi}.
\end{proof}

Combining Lemma~\ref{nondecreasing} with Example~\ref{graphs}, we obtain the principal result
of this section.

\begin{prop}
\label{nonneg}
Assume that there exists a non-negative Legendrian isotopy connecting the zero
section of $\Jet^1(L)$ with the graph $\Lambda^f$ of the $1$-jet of a smooth
function~$f$. Then $f\ge 0$ everywhere on $L$.\qed
\end{prop}

Here are a couple of immediate corollaries of Proposition~\ref{nonneg}.

\begin{cor}
\label{nonneg1}
If there exists a non-negative Legendrian isotopy connecting $\Lambda^h$
to $\Lambda^f$, then $h\le f$ everywhere on~$L$.
\end{cor}

\begin{proof}
The `downshift' map $T^h(q,p,u):=(q,p-\frac{\partial h}{\partial q}, u-h)$ is a contactomorphism of $\Jet^1(L)$.
Clearly, $T^h(\Lambda^h)=\Lambda^0$ and $T^h(\Lambda^f)=\Lambda^{f-h}$. Thus,
$T^h$ will take a non-negative Legendrian isotopy between $\Lambda^h$
and $\Lambda^f$ to a non-negative Legendrian isotopy between
$\Lambda^0$ and $\Lambda^{f-h}$. Hence, $f-h\ge 0$ by Proposition~\ref{nonneg}.
\end{proof}

\begin{cor}
\label{nonneg2}
A non-negative Legendrian isotopy connecting the zero section with
itself is constant.
\end{cor}

\begin{proof}
For any non-constant Legendrian isotopy  $\{\Lambda_t\}_{t\in[0,1]}$
of the zero section, there exists $t'>0$ such that $\Lambda_{t'}=\Lambda^h$ for
a function~$h$ not identically equal to zero. (Note that $\Lambda_t$ is the
graph of the $1$-jet of a function for all sufficiently small $t\ge 0$.)
Assume that the isotopy is non-negative.
Proposition~\ref{nonneg} shows that $h\ge 0$. On the other hand, applying
Corollary~\ref{nonneg1} to the isotopy $\{\Lambda_t\}_{t\in[t',1]}$, we see that $h\le 0$.
So $h\equiv 0$, a contradiction.
\end{proof}

\subsection{Application of the hodograph transformation}
\label{hodograph}
Let $\langle\cdot,\cdot\rangle$ denote the standard scalar product on $\R^m$
and let $S^{m-1}\subset\R^m$ be the unit sphere.
The map
$$
\R^m\times S^{m-1}\ni (x,q)\longmapsto\langle q,\cdot \rangle \in ST_x^*\R^m
$$
provides a trivialisation of $ST^*\R^m$.
The {\it hodograph transformation\/} is then defined by the formula
\begin{equation}
\label{hodo}
ST^*\R^m\ni(x,q) \longmapsto
(q,\langle x,\cdot \rangle|_{T_qS^{m-1}}, \langle x,q\rangle)\in\Jet^1(S^{m-1}).
\end{equation}
It is easy to see that this map is a contactomorphism of the standard contact structures on $ST^*\R^m$ and $\Jet^1(S^{m-1})$
(see~\cite[pp.\,48--49]{Ar}).

In the case $m=2$, we can trivialise $\Jet^1(S^1)$ using the angle coordinate~$\phi$ on $S^1$
and the corresponding momentum coordinate on the fibre of $T^*S^1$.
Formula~(\ref{hodo}) then becomes
\begin{equation}
\label{hodo2}
(x_1,x_2,\phi)\longmapsto (\phi, -x_1\sin\phi+x_2\cos\phi, x_1\cos\phi +x_2\sin\phi).
\end{equation}

It should be clear from formula~(\ref{hodo}) that the hodograph image
of the fibre of $ST^*\R^m$ over a point $x\in\R^m$ is the graph of
the $1$-jet of the function $q\mapsto\langle x,q\rangle$ on $S^{m-1}$.
In particular, the fibre over the origin is mapped to the zero section
of $\Jet^1(S^{m-1})$.

\begin{cor}
\label{fibres}
If a non-negative Legendrian isotopy connects two fibres of $ST^*\R^m$,
then the fibres coincide and the isotopy is constant.
\end{cor}

\begin{proof}
Applying a parallel shift in $\R^m$, we may assume that the first fibre is the fibre over the origin.
Then the image of our isotopy under the hodograph transformation is a non-negative
Legendrian isotopy in $\Jet^1(S^{m-1})$ connecting the zero section $\Lambda^0$ to the graph
$\Lambda^f$ of the $1$-jet of the function $f(q)=\langle x,q\rangle$ for some $x\in\R^m$.
Proposition~\ref{nonneg} shows that $f\ge 0$, which is only possible if $x=0$ because
otherwise $f(-x/|x|)=-|x|<0$. Thus, $\Lambda^f=\Lambda^0$ and the isotopy is constant by
Corollary~\ref{nonneg2}.
\end{proof}

\begin{cor}
\label{fibresdown}
Let $M$ be a manifold smoothly covered by an open subset $\wt M\subset\R^m$.
If a non-negative Legendrian isotopy connects two fibres of $ST^*M$,
then the fibres coincide and the isotopy is constant.
\end{cor}

\begin{proof}
Let $ST^*\wt M\to ST^*M$ be the fibrewise covering associated to the covering $\wt M\to M$.
A non-negative Legendrian isotopy connecting two fibres of $ST^*M$ lifts to
a non-negative Legendrian isotopy connecting two fibres of  $ST^*\wt M\subset ST^*\R^m$.
Thus, the result follows from Corollary~\ref{fibres}.
\end{proof}

\begin{rem}
\label{exhaust}
Somewhat more generally, one can assume that every compact subset of $M$
is contained in an open set smoothly covered by an open subset of~$\R^m$.
\end{rem}

\begin{exm}
There exist manifolds $M$ for which the assertion of Corollary~\ref{fibresdown} is false.
The simplest example is provided by the sphere $S^m$. Indeed, moving a fibre of $ST^*S^m$
along the (co-)geodesic flow of the standard round metric $\ov g$ defines a positive Legendrian isotopy
connecting it to the fibre over the antipodal point and then to itself.
Note that a very similar argument shows that the Legendrian Low conjecture is false
for the product Lorentz manifold $(S^{m}\times \R, \ov g\oplus -dt^2)$
with Cauchy surface $S^m$, see~\cite[Example 3]{ChernovRudyak}.
\end{exm}

\subsection{Proof of the Legendrian Low conjecture}
\label{proofleg}
Let $x$ and $y$ be causally related points in a globally hyperbolic spacetime with a
Cauchy surface $M$ covered by an open subset of~$\R^m$. Identify
the space of future pointing null geodesics $\mathfrak N$
with~$ST^*M$. Since intersecting skies are linked by definition,
we may assume that $x$ and $y$ do not lie on the same null geodesic.
We may also assume that $y$ lies in the causal past of $x$.
By Proposition~\ref{causal} there is a {\it non-negative\/}
Legendrian isotopy $\{\Lambda_t\}_{t\in[0,1]}$ in $ST^*M$
such that $\Lambda_0=\mathfrak S_x$ and $\Lambda_1=\mathfrak S_y$.

Suppose that $\mathfrak S_x$ and $\mathfrak S_y$ are Legendrian unlinked, i.\,e.,
that the link $\mathfrak S_x\sqcup\mathfrak S_y$ is Legendrian isotopic to
the link $F\sqcup F'$ formed by  two different fibres of $ST^*M$.
By the Legendrian isotopy extension theorem,
there exists a contactomorphism $\Psi:ST^*M\to ST^*M$ such that
$\Psi(\mathfrak S_x\sqcup\mathfrak S_y)=F\sqcup F'$.
Hence, $\{\Psi(\Lambda_t)\}_{t\in[0,1]}$ is a non-negative Legendrian isotopy
connecting two different fibres of $ST^*M$, which contradicts Corollary~\ref{fibresdown}.
Thus, $\mathfrak S_x$ and $\mathfrak S_y$ are Legendrian linked.\qed

\subsection{Proof of Theorem~\ref{thmc}}
\label{proofthmc}
Consider again a globally hyperbolic spacetime with a Cauchy surface $M$
covered by an open subset of $\R^m$. Identify $\mathfrak N$ with~$ST^*M$.

Let $x$ and $y$ be causally related points with disjoint skies.
Suppose that the links $\mathfrak S_x\sqcup\mathfrak S_y$
and $\mathfrak S_y\sqcup\mathfrak S_x$ are Legendrian isotopic.
By the Legendrian isotopy extension theorem, there exists a contactomorphism $\Psi$
such that $\Psi(\mathfrak S_x\sqcup\mathfrak S_y)=\mathfrak S_y\sqcup\mathfrak S_x$.
Assume that $y$ lies in the causal past of $x$ (otherwise rename the points).
Let $\{\Lambda_t\}_{t\in[0,1]}$ be a non-negative Legendrian isotopy
in $ST^*M$ connecting $\mathfrak S_x$ to $\mathfrak S_y$
provided by Proposition~\ref{causal}. Then $\{\Psi(\Lambda_t)\}_{t\in[0,1]}$ is
a non-negative Legendrian isotopy connecting $\mathfrak S_y$ to~$\mathfrak S_x$.
Composing these two isotopies, we obtain a non-constant non-negative
Legendrian isotopy connecting $\mathfrak S_x$ to itself.

Recall now that any sky $\mathfrak S_x$ is Legendrian isotopic to a fibre of $ST^*M$.
By the Legendrian isotopy extension theorem, there exists a contactomorphism $\Phi$
taking $\mathfrak S_x$ to that fibre. The Legendrian isotopy connecting $\mathfrak S_x$ to itself
constructed above is taken by $\Phi$ to a non-constant non-negative Legendrian isotopy connecting
the fibre to itself, which is impossible by Corollary~\ref{fibresdown}.
Hence, the links
$\mathfrak S_x\sqcup\mathfrak S_y$ and $\mathfrak S_y\sqcup\mathfrak S_x$
cannot be Legendrian isotopic.\qed

\begin{rem}
\label{generalised}
The proof shows that the links $\mathfrak S_x\sqcup\mathfrak S_y$ and $\mathfrak S_y\sqcup\mathfrak S_x$
are not even contactomorphic in $\mathfrak N$. Similarly, the proof of Theorem~\ref{thma} in \S\ref{proofleg}
shows that the link formed by the skies of causally related points $x,y\in X$ is not contactomorphic
to any link formed by the skies of a pair of causally unrelated points.
\end{rem}

\begin{rem}
\label{extended}
The proofs of Theorems~\ref{thma} and~\ref{thmc} in~\S\S\ref{proofleg}--\ref{proofthmc}
work for any globally hyperbolic spacetime such that the assertion of Corollary~\ref{fibresdown}
is true for its Cauchy surface~$M$.
\end{rem}

\subsection{Proof of the Low conjecture for (2+1)-dimensional spacetimes}
\label{prooftop}
Consider the (nonoriented) link in $ST^*\R^2$ formed by two distinct fibres.
In the terminology of~\cite{DG}, the image of this link under
the hodograph transformation~(\ref{hodo}) is the $(1,1)$-cable link in $\Jet^1(S^1)$.
(Indeed, it is clear from formula~(\ref{hodo2}) that if the image of the first fibre
is the zero section, then the image of the second one goes along
it once and makes a single turn around it in the $(p,u)$-plane.)
According to the main result of~\cite{DG}, Legendrian links smoothly
isotopic to the $(1,1)$-cable link are classified up to Legendrian
isotopy by the classical invariants of their components.

Now let  $\ss$ be a globally hyperbolic $(2+1)$-dimensional spacetime
with Cauchy surface~$M$ covered by $\R^2$.
Let $\mathfrak S_x\sqcup\mathfrak S_y$ be a link in $\mathfrak N\cong ST^*M$
formed by the skies of two causally related points. We may assume that $y\in J^-(x)$.
Let $\{\Lambda_t\}_{t\in[0,1]}$ be a non-negative Legendrian isotopy in $ST^*M$
connecting $\mathfrak S_x$ to $\mathfrak S_y$ provided by Proposition~\ref{causal}.
Consider the covering $ST^*\R^2\to ST^*M$ associated to the covering of~$M$ by~$\R^2$.
Lift the isotopy $\{\Lambda_t\}_{t\in[0,1]}$ to a non-negative Legendrian isotopy
$\{\wt\Lambda_t\}_{t\in[0,1]}$ in $ST^*\R^2$ and set $\wt{\mathfrak S}_x=\wt\Lambda_0$
and $\wt{\mathfrak S}_y=\wt\Lambda_1$.

Now we argue by contradiction. Suppose that the link $\mathfrak S_x\sqcup\mathfrak S_y$
is smoothly isotopic to a pair of fibres in $ST^*M$. Since this isotopy lifts to $ST^*\R^2$,
it follows that $\wt{\mathfrak S}_x\sqcup\wt{\mathfrak S}_y$ is smoothly isotopic
to a pair of fibres in~$ST^*\R^2$. Each component of $\wt{\mathfrak S}_x\sqcup\wt{\mathfrak S}_y$
is Legendrian isotopic to the fibre of $ST^*\R^2$ and therefore has the same classical invariants (in $\Jet^1(S^1)$).
Hence, $\wt{\mathfrak S}_x\sqcup\wt{\mathfrak S}_y$ is Legendrian isotopic to
a pair of fibres by the aforementioned result from~\cite{DG}.
Applying the Legendrian isotopy extension theorem in the same way
as in~\S\ref{proofleg}, we see that this contradicts Corollary~\ref{fibres}.\qed

\begin{rem}
Let $\wt \ss\to \ss$ be the covering corresponding to the covering $\R^2\to M$.
Then $(\wt\ss, \wt g)$ with the induced Lorentz metric is a globally hyperbolic spacetime
with a spacelike Cauchy surface~$\wt M\cong\R^2$, see~\cite[Proof of Theorem~14]{ChernovRudyak}.
For a pair of causally related points $x,y\in \ss$ with $y\in J^-(x)$, we can pick a pair of lifts
$\wt x, \wt y\in \wt \ss$ connected by the lift of a past pointing curve connecting $x$ to~$y$.
In this way, Theorem~\ref{thmb} can be reduced to the case when the Cauchy
surface is itself diffeomorphic to $\R^2$. A similar argument may be used to reduce
Theorem~\ref{thma} to the case when the Cauchy surface is diffeomorphic to an open subset of~$\R^m$.
Note also that the Legendrian circles $\wt{\mathfrak S}_x, \wt{\mathfrak S}_y\subset ST^*\R^2$ appearing
in the proof of Theorem~\ref{thmb} correspond to the skies ${\mathfrak S}_{\wt x}, {\mathfrak S}_{\wt y}$.
\end{rem}

\subsection{Application to relativity}
\label{scifi}
Suppose that $M$ is a manifold smoothly covered by an open subset of $\R^m$.
It follows from Theorems~\ref{thma} and~\ref{thmc} that it is
impossible to find two globally hyperbolic Lorentz metrics $g_1$ and $g_2$
on $\ss=M\times \R$ and two pairs of events $x_1, y_1\in (\ss, g_1)$
and $x_2, y_2\in (\ss, g_2)$ such that
\begin{itemize}
\item[a)] $M_0=M\times 0$ is a spacelike Cauchy surface for both $(\ss, g_1)$ and ($\ss, g_2)$;
\item[b)] the Legendrian links in $ST^*M_0$ corresponding to $(x_1, y_1)$ and $(x_2,y_2)$ are the same
(in any natural sense, see Remark~\ref{generalised});
\item[c)] $x_1, y_1$ are causally related and $x_2, y_2$ are causally unrelated
(or $y_1\in J^+(x_1)$ and $y_2\not \in J^+(x_2))$.
\end{itemize}

The Legendrian links in (b) are completely determined by the intersections
of the exponents of the null cones at $x_1,y_1$ and $x_2,y_2$ with an arbitrarily thin
(in fact, infinitesimal) neighbourhood of the Cauchy surface~$M_0$. Hence, the result
admits the following physical or, perhaps, `science fictional' interpretation.
Even if you are given absolute freedom to change the past, you will not be able
to destroy causal relations between events so that the change cannot be observed on spacelike
Cauchy surfaces in the immediate present.

\smallskip
\noindent
{\bf Acknowledgment.} The authors are grateful to Hansj\"org Geiges,
Robert Low, and the anonymous referee for useful comments.


\begin{thebibliography}{99}
\bibitem{AgolTalk}
I.~Agol, {\it The geometrization conjecture and universal covers of $3$-manifolds}, a talk at the Cornell
2004 Topology Festival. The transparencies of the talk are available at
{\tt http://www2.math.uic.edu/\verb'~'agol/cover/cover01.html}
\bibitem{Ar}
V. I. Arnold, {\it Topological invariants of plane curves and caustics},
University Lecture Series, 5. American Mathematical Society, Providence, RI, 1994.
\bibitem{ArnoldProblem}
V.~I.~Arnold, {\em Problems,\/} written down by S.~Duzhin, September 1998,
available electronically at {\tt http://www.pdmi.ras.ru/\verb'~'arnsem/Arnold/prob9809.ps.gz}
\bibitem{ArnoldProblemBook}
V.~I.~Arnold, {\em Arnold's problems}, Translated and revised edition of the 2000 Russian original.
Springer-Verlag, Berlin; PHASIS, Moscow, 2004.
\bibitem{BeemEhrlichEasley}
J.~K.~Beem, P.~E.~Ehrlich, K.~L.~Easley, {\em Global Lorentzian geometry},
Second edition. Monographs and Textbooks in Pure and Applied Mathematics {\bf 202},
Marcel Dekker, Inc., New York, 1996.
\bibitem{BerardBergery}
L.~B\'erard-Bergery,
{\em Quelques exemples de vari\'et\'es riemanniennes o\`u toutes les
g\'eod\'esiques issues d'un point sont ferm\'ees et de m\^ eme longueur,
suivis de quelques r\'esultats sur leur topologie},
Ann. Inst. Fourier (Grenoble) {\bf 27} (1977), 231--249.
\bibitem{BernalSanchez}
A.~Bernal, M.~Sanchez, {\em On smooth Cauchy hypersurfaces and Geroch's splitting theorem},
Comm.~Math.~Phys.~{\bf 243} (2003), 461--470.
\bibitem{BernalSanchezMetricSplitting}
A.~Bernal, M.~Sanchez, {\em Smoothness of time functions and the
metric splitting of globally hyperbolic space-times},
Comm.~Math.~Phys.~{\bf 257} (2005), 43--50.
\bibitem{BernalSanchezFurther}
A.~Bernal, M.~Sanchez, {\em Further results on the smoothability
of Cauchy hypersurfaces and Cauchy time functions},
Lett.~Math.~Phys.~{\bf 77} (2006), 183--197.
\bibitem{BernalSanchezCausal}
A.~Bernal, M.~Sanchez, {\it Globally hyperbolic spacetimes can be defined as ``causal'' instead of
``strongly causal''},  Class.~Quant.~Grav.~{\bf 24} (2007), 745--750.
\bibitem{Besse}
A.~L.~Besse, {\em Manifolds all of whose geodesics are closed},
with appendices by D.~B.~A.~Epstein, J.-P.~Bourguignon,
L.~B\'erard-Bergery, M.~Berger and J.~L.~Kazdan.
Ergebnisse der Mathematik und ihrer Grenzgebiete, 93. Springer-Verlag, Berlin-New York, 1978.
\bibitem{Bh} M. Bhupal, {\it A partial order on the group of contactomorphisms of
$\mathbb R\sp {2n+1}$ via generating functions}, Turkish J. Math. {\bf 25} (2001), 125--135.
\bibitem{Cha} M. Chaperon, {\it On generating families},  The Floer memorial volume,  283--296,
Progr. Math., 133, Birkh\"auser, Basel, 1995.
\bibitem{Che} Yu. V. Chekanov, {\it Critical points of quasifunctions, and generating families of Legendrian
manifolds}, Funktsional. Anal. i Prilozhen.  {\bf 30}:2 (1996), 56--69 (Russian);
English transl. in Funct. Anal. Appl. {\bf 30}:2 (1996), 118--128.
\bibitem{ChernovRudyak}
V.~Chernov (Tchernov), Yu.~Rudyak,
{\it Linking and causality in globally hyperbolic space-times},
Comm. Math. Phys. {\bf 279} (2008), 309--354.
\bibitem{CFP}
V. Colin, E. Ferrand, P. Pushkar,
{\it Positive isotopies of Legendrian submanifolds},
Preprint available at {\tt http://people.math.jussieu.fr/\verb'~'ferrand/publi/PIL.pdf}
\bibitem{DG}
F. Ding, H. Geiges, {\it Legendrian helix and cable links}, Preprint {\tt arXiv:math/0611080v2},
Commun. Contemp. Math., to appear.
\bibitem{ElKiPo}
Y. Eliashberg, S. S. Kim, L. Polterovich, {\it Geometry of contact transformations and domains\/{\rm :}
orderability versus squeezing}, Geom. Topol. {\bf 10} (2006), 1635--1747.
\bibitem{ElPo}
Y. Eliashberg, L. Polterovich, {\it Partially ordered groups and geometry of contact transformations},
Geom. Funct. Anal. {\bf 10} (2000), 1448--1476.
\bibitem{Ge}
H. Geiges, {\it An introduction to contact topology},
Cambridge Studies in Advanced Mathematics, 109. Cambridge University Press, Cambridge, 2008.
\bibitem{Geroch}
R.~P.~Geroch, {\em Domain of dependence}, J.~Math.~Phys. {\bf 11} (1970), 437--449.
\bibitem{HawkingEllis}
S.~W.~Hawking, G.~F.~R.~Ellis, {\em The large scale structure of space-time},
Cambridge Monographs on Mathematical Physics, No.~1,
Cambridge University Press, London-New York, 1973.
\bibitem{Low0}
R.~J.~Low, {\em Causal relations and spaces of null geodesics},
PhD Thesis, Oxford University (1988).
\bibitem{Low1}
R.~J.~Low, {\em Twistor linking and causal relations},
Classical Quantum Gravity {\bf 7} (1990), 177--187.
\bibitem{Low3}
R.~J.~Low, {\em Twistor linking and causal relations in exterior Schwarzschild space},
Classical Quantum Gravity {\bf 11} (1994), 453--456.
\bibitem{LowLegendrian}
R.~J.~Low, {\em Stable singularities of wave-fronts in general relativity},
J.~Math.~Phys.~{\bf 39} (1998), 3332--3335.
\bibitem{LowNullgeodesics}
R.~J.~Low, {\em The space of null geodesics}, Proceedings of the
Third World Congress of Nonlinear Analysts, Part 5 (Catania, 2000).
Nonlinear Anal. {\bf 47} (2001), 3005--3017.
\bibitem{LowRefocussing}
R.~J.~Low, {\em The space of null geodesics {\rm(}and a new causal boundary\/{\rm )}},
Lecture Notes in Physics {\bf 692}, Springer, Berlin Heidelberg New York, 2006, pp.\,35--50.
\bibitem{NatarioTod}
J.~Nat\'ario, P.~Tod,
{\em Linking, Legendrian linking and causality},
Proc.~London Math.~Soc.~(3) {\bf 88} (2004), 251--272.
\bibitem{ONeill}
B.~O'Neill, {\em Semi-Riemannian geometry. With applications to
relativity}, Pure and Applied Mathematics, 103. Academic Press,
Inc. [Harcourt Brace Jovanovich, Publishers], New York, 1983.
\bibitem{Perelman1}
G.~Perelman, {\em The entropy formula for the Ricci flow and its
geometric applications}, Preprint {\tt math.DG/0211159}.
\bibitem{Perelman2}
G.~Perelman, {\em Ricci flow with surgery on three-manifolds},
Preprint {\tt math.DG/0303109}.
\bibitem{Tr} L. Traynor, {\it Legendrian circular helix links},
Math. Proc. Cambridge Philos. Soc. {\bf 122} (1997), 301--314.
\bibitem{Vi} C. Viterbo, {\it Symplectic topology as the geometry of generating functions},
Math. Ann. {\bf 292} (1992), 685--710.
\end{thebibliography}
\end{document}